\documentclass[11pt,a4paper]{amsart}
\usepackage{graphicx}
\usepackage{amssymb}
\usepackage{amsmath}
\usepackage{amsthm}
\usepackage{amscd}
\usepackage[all,2cell]{xy}
 \UseAllTwocells \SilentMatrices
\newtheorem{thm}{Theorem}[section]

\newtheorem{cor}[thm]{Corollary}
\newtheorem{lem}[thm]{Lemma}
\newtheorem{exm}[thm]{Example}

\newtheorem{prop}[thm]{Proposition}
\theoremstyle{definition}

\theoremstyle{remark}
\newtheorem{rem}[thm]{\bf Remark}
\numberwithin{equation}{section}

\begin{document}
\title[Generalized Serre duality]{Generalized Serre duality}
\author[X. W. Chen
] {Xiao-Wu Chen}
\thanks{This project was supported by Alexander von Humboldt Stiftung and National Natural Science
Foundation of China (No. 10971206).}
\thanks{E-mail: xwchen$\symbol{64}$mail.ustc.edu.cn}
\date{\today}
\maketitle
\dedicatory{}%
\commby{}%


\begin{abstract}
We introduce a notion of generalized Serre duality on a Hom-finite
Krull-Schmidt triangulated category $\mathcal{T}$. This duality
induces the generalized Serre functor on $\mathcal{T}$, which is a
linear triangle equivalence between two thick triangulated
subcategories of $\mathcal{T}$. Moreover, the  domain of the
generalized Serre functor  is the smallest additive subcategory of
$\mathcal{T}$ containing all the indecomposable objects  which
appear as the third term of an Auslander-Reiten triangle in
$\mathcal{T}$; dually, the range of the generalized Serre functor is
the smallest additive subcategory of $\mathcal{T}$ containing all
the indecomposable objects  which appear as the first term of an
Auslander-Reiten triangle in $\mathcal{T}$.

   We compute explicitly the generalized Serre duality on the
bounded derived categories of artin algebras  and of certain
noncommutative projective schemes in the sense of Artin and Zhang.
We obtain a characterization of Gorenstein algebras: an artin
algebra $A$ is Gorenstein if and only if the bounded homotopy
category  of finitely generated projective $A$-modules has Serre
duality in the sense of Bondal and Kapranov.
\end{abstract}

\section{Introduction}

Throughout $R$ will be a commutative artinian ring. Let
$\mathcal{T}$ be an $R$-linear triangulated category which is
\emph{Hom-finite} and \emph{Krull-Schmidt}. Here, an $R$-linear
category is Hom-finite provided that all its Hom spaces are finite
generated $R$-modules; an additive category is Krull-Schmidt
provided that each object is a finite sum of indecomposable objects
with local endomorphism rings.

The Auslander-Reiten theory for a Hom-finite Krull-Schmidt
triangulated category $\mathcal{T}$ was initiated by Happel
(\cite{H1, H2}) which now plays an important role in the
representation theory of artin algebras. The central notion is
\emph{Auslander-Reiten triangle}, which by definition is a
distinguished triangle $X\stackrel{u}\rightarrow Y
\stackrel{v}\rightarrow Z \rightarrow X[1]$ with $X, Z$
indecomposable such that $u$ (\emph{resp}. $v$) is left
(\emph{resp}. right) almost-split. Here $[1]$ denotes the
translation functor on $\mathcal{T}$. For an indecomposable object
$X$ there exists at most one, up to isomorphism, Auslander-Reiten
triangle which contains $X$ as the first term; dually, for an
indecomposable object $Z$ there exists at most one Auslander-Reiten
triangle which contains $Z$ as the third term. The triangulated
category $\mathcal{T}$ is said to \emph{have enough Auslander-Reiten
triangles} provided that for each indecomposable object $X$ there
exist  two Auslander-Reiten triangles which contain $X$ as the first
term and the third term, respectively. A fundamental result due to
Reiten and Van den Bergh states that $\mathcal{T}$ has enough
Auslander-Reiten triangles if and only if it has Serre duality
(\cite[Theorem I.2.4]{RV}). Here the Serre duality is in the sense
of Bondal and Kapranov (\cite{BK}).

Observe that in general a Hom-finite Krull-Schmidt triangulated
category $\mathcal{T}$ will not have enough Auslander-Reiten
triangles. For example, Happel proved that the bounded derived
category of an artin algebra has enough Auslander-Reiten triangles
if and only if the algebra has finite global dimension (\cite{H2}).
In this situation, the category $\mathcal{T}$ will not have Serre
duality.

In the present paper, we introduce a notion of  \emph{generalized
Serre duality} on a Hom-finite Krull-Schmidt triangulated category
$\mathcal{T}$. This generalized duality induces the
\emph{generalized Serre functor} of $\mathcal{T}$, which is a linear
triangle equivalence between  certain \emph{thick} triangulated
subcategories of $\mathcal{T}$. Here, thick subcategories mean
subcategories which are closed under direct summands. We put
Reiten-Van den Bergh's theorem in this general setting. It turns out
that this generalized Serre duality is closely related the
Auslander-Reiten triangles in $\mathcal{T}$.

The notion of generalized Serre duality applies to an arbitrary
Hom-finite $R$-linear category. We need some notation. Let
$\mathcal{C}$ be a Hom-finite $R$-linear category. For two objects
$X$ and $Y$, we write $(X, Y)$ for the Hom space ${\rm
Hom}_\mathcal{C}(X, Y)$. Denote by $(X, -)$ (\emph{resp}. $(-, X)$)
the representable functor ${\rm Hom}_\mathcal{C}(X, -)$
(\emph{resp}. ${\rm Hom}_\mathcal{C}(-, X)$ ) from $\mathcal{C}$ to
the category $R\mbox{-mod}$ of finitely generated  $R$-modules.
Denote by $D$ the Matlis duality on $R\mbox{-mod}$. Recall that
$D={\rm Hom}_R(-, E)$ where $E$ is the minimal injective cogenerator
for $R$; see \cite[Chapter II.3]{ARS}. By $D(X, -)$ (\emph{resp}.
$D(-, X)$) we mean the composite functor $D{\rm Hom}_\mathcal{C}(X,
-)$ (\emph{resp}. $D{\rm Hom}_\mathcal{C}(-, X)$).

We define two full subcategories of $\mathcal{C}$ as follows
\begin{align*}
\mathcal{C}_r=\{X \in \mathcal{C}\; |\; D(X, -) \mbox{ is
representable}  \} \mbox{ and } \mathcal{C}_l=\{X \in
\mathcal{C}\;|\; D(-, X) \mbox{ is representable}\}.
\end{align*}
It turns out that there exists a unique functor  $S\colon
\mathcal{C}_r \rightarrow \mathcal{C}_l$ such that there is a
$R$-linear isomorphism
$$\phi_{X, Y}\colon D(X, Y) \simeq (Y, S(X))$$
for each $X \in \mathcal{C}_r$ and $Y \in \mathcal{C}$; the
isomorphism $\phi_{X, Y}$ is required to be natural both in $X$ and
$Y$. In fact, the functor $S$ is an $R$-linear equivalence of
categories; see Appendix A.1. We will call the functor $S$ the
\emph{generalized Serre functor} on $\mathcal{C}$, and we refer to
$\mathcal{C}_r$ (\emph{resp}. $\mathcal{C}_l$) as the \emph{domain}
(\emph{resp}. the \emph{range}) of the generalized Serre functor. We
will call this set of data the \emph{generalized Serre duality} on
$\mathcal{C}$. Observe that the category $\mathcal{C}$ has Serre
duality in the sense of \cite{BK} if and only if
$\mathcal{C}_r=\mathcal{C}=\mathcal{C}_l$. For details, see
Appendix.

 \par \vskip 5pt

Here is our main theorem.   \\

\noindent {\bf Main Theorem.}\quad  \emph{Let $\mathcal{T}$ be a
Hom-finite Krull-Schmidt triangulated category. Denote by $S \colon
\mathcal{T}_r \rightarrow \mathcal{T}_l$ its generalized Serre
functor. Then both $\mathcal{T}_r$ and $\mathcal{T}_l$ are thick
triangulated subcategories of $\mathcal{T}$. Moreover, we have
\begin{enumerate}\item there is a natural isomorphism $\eta_X \colon
S(X[1])\rightarrow S(X)[1]$ for each $X \in \mathcal{T}_r$ such that
the pair $(S, \eta)$ is a triangle equivalence between
$\mathcal{T}_r$ and $\mathcal{T}_l$;
\item  an indecomposable object $X$ in $\mathcal{T}$ belongs to
$\mathcal{T}_l$ (resp. $\mathcal{T}_r$) if and only if there is an
Auslander-Reiten triangle in $\mathcal{T}$ containing $X$ as the
first term (resp. the third term).\end{enumerate}}

\vskip 10pt

Let us remark that the statement (1) is a generalization of a result
due to Bondal and Kapranov (\cite{BK}; also see \cite{Bo,Huy}). The
main theorem is a combination of Propositions
\ref{prop:domainrange}, \ref{prop:2.7}, \ref{prop:2.8} and
\ref{prop:A.1}.

\par \vskip 5pt

The paper is organized as follows. In section 2, we divide the proof
of the main theorem into proving several propositions. Section 3 is
devoted to computing the generalized Serre duality for the bounded
derived categories of certain abelian categories explicitly; see
Theorem \ref{thm:3.5} and Theorem \ref{thm:3.10}. The computational
results are based on results by Happel (\cite{H2}), and by de
Naeghel and Van den Bergh ({\cite{NV}). We obtain, as a byproduct, a
seemingly new characterization of Gorenstein algebras: an artin
algebra $A$ is Gorenstein if and only if the bounded homotopy
category  of finitely generated projective $A$-modules has Serre
duality;  see Corollary \ref{cor:3.9}. The appendix explains the
basic notions and results on generalized Serre duality on a
Hom-finite linear category.
\par

 For triangulated categories, we refer to \cite{H1, Har, Ne}. For the notion of
Auslander-Reiten triangles, we refer to \cite{H1, H2}.
\par \vskip 10pt

\section{Proof of Main Theorem}

In this section, we divide the proof of the main theorem into
proving several propositions. We make preparation on a study of
coherent functors on triangulated categories. We collect the
relevant results  which seems to be scattered in the literature.

\vskip 5pt

 Let $\mathcal{C}$ be an additive
category. Denote by $(\mathcal{C}^{\rm op}, {\rm Ab})$ the large
category of contravariant additive  functors from $\mathcal{C}$ to
the category  Ab of abelian groups. Note that although
$(\mathcal{C}^{\rm op}, {\rm Ab})$ is usually not a category, but it
still makes sense to mention the exact sequences of functors, in
particular, the notions of kernel, cokernel and extension in
$(\mathcal{C}^{\rm op}, {\rm Ab})$ make sense. Hence sometimes we
even pretend that $(\mathcal{C}^{\rm op}, {\rm Ab})$ is an abelian
category.

Recall that the Yoneda embedding functor ${\bf p}\colon \mathcal{C}
\rightarrow (\mathcal{C}^{\rm op}, {\rm Ab})$ is defined by ${\bf
p}(C)= (-, C)$, where $(-, C)$ denotes the representable functor
${\rm Hom}_{\mathcal{C}}(-, C)$. Then Yoneda  Lemma implies that
${\bf p}$ is fully faithful; moreover, representable functors are
projective objects in $(\mathcal{C}^{\rm op}, {\rm Ab})$. Recall
that a functor $F\in (\mathcal{C}^{\rm op}, {\rm Ab}) $ is said to
be \emph{coherent} provided that there exists an exact sequence of
functors
$$(-, C_0)
\longrightarrow (-, C_1)\longrightarrow F \longrightarrow 0,$$ where
$C_0, C_1 \in \mathcal{C}$. Note that such a sequence can be viewed
as a projective presentation of $F$ in $(\mathcal{C}^{\rm op}, {\rm
Ab})$. Denote by $\widehat{\mathcal{C}}$ the full subcategory of
$(\mathcal{C}^{\rm op}, {\rm Ab})$ consisting of coherent functors
on $\mathcal{C}$. Note that the category $\widehat{\mathcal{C}}$ of
coherent functors is an additive category with small Hom sets.
\par

Recall that a \emph{pseudokernel} of a morphism $f\colon
C_0\rightarrow C_1$ in $\mathcal{C}$ is a morphism $k\colon K
\rightarrow C_0$ such that $f\circ k=0$ and that each morphism
$k'\colon K \rightarrow C_0$ with $f\circ k'=0$ factors through $k$,
in other words, the morphism $k$ makes the sequence $(-, K)
\stackrel{(-, k)}{\rightarrow} (-,C_0) \stackrel{(-,
f)}{\rightarrow} (-,C_1)$ of functors exact. We say that the
category $\mathcal{C}$ \emph{has pseudokernels} provided that each
morphism has a pseudokernel. Dually, one defines
\emph{pseudocokernels}. \par \vskip 5pt

Let us recall a basic result on coherent functors (\cite{Au1, Au2};
also see \cite{Kr}).

 \begin{lem}\label{lem:coherent}
Let $\mathcal{C}$ be an additive category. Then the category
$\widehat{\mathcal{C}}$ of coherent functors is closed under
cokernels and extensions in $(\mathcal{C}^{\rm op}, {\rm Ab})$. In
particular, $\widehat{\mathcal{C}}$ is closed under taking direct
summands. Moreover, $\widehat{\mathcal{C}} $ is an abelian
subcategory of $(\mathcal{C}^{\rm op}, {\rm Ab})$ if and only if
$\mathcal{C}$ has pseudokernels. \hfill $\square$
\end{lem}

The following is well known, as pointed out in \cite[Example
4.1(2)]{Kr}.

\begin{lem}\label{lem:trianglecoherent}
Let $\mathcal{T}$ be a triangulated category. Then $\mathcal{T}$ has
pseudokernels and pseudocokernels.
\end{lem}

\begin{proof} Given a morphism $v: Y \rightarrow Z$ in $\mathcal{T}$, we
take a distinguished triangle
$$X\stackrel{u}{\longrightarrow} Y \stackrel{v}{\longrightarrow} Z
\stackrel{w}{\longrightarrow} X[1].$$
 Since for any object $C\in
\mathcal{T}$, the functors $(C, -)$ and $(-, C)$ are cohomological
(see \cite[p.23]{Har}), it is easily verified  that $u$ and $w$ are
a pseudokernel and a pseudocokernel of $v$, respectively.
\end{proof}

In what follows $\mathcal{T}$ will  be a triangulated category. We
will introduce a {\em two-sided resolution} for any coherent functor
on $\mathcal{T}$. Given a functor $F \in \widehat{\mathcal{T}}$, we
have a presentation
\begin{align*}
(-, Y)\stackrel{(-, v)}{\longrightarrow} (-, Z) \longrightarrow F
\longrightarrow 0.
\end{align*}
Take a distinguished triangle $X\stackrel{u}{\rightarrow} Y
\stackrel{v}{\rightarrow} Z \stackrel{w}{\rightarrow} X[1].$ So we
have a long exact sequence of functors
$$ \cdots
\longrightarrow (-, Z[-1]) \stackrel{(-, w[-1])}{\longrightarrow}
(-, X) \stackrel{(-, u)}{\longrightarrow} (-, Y) \stackrel{(-,
v)}{\longrightarrow} (-, Z) \stackrel{(-, w)}{\longrightarrow} (-,
X[1]) \longrightarrow \cdots $$

Note that $F \simeq {\rm Coker}(-, v)$. Hence we have the
following two exact sequences
$$ \cdots
\longrightarrow (-, Z[-1]) \stackrel{(-, w[-1])}{\longrightarrow}
(-, X) \stackrel{(-, u)}{\longrightarrow} (-, Y) \stackrel{(-,
v)}{\longrightarrow} (-, Z) \longrightarrow F \longrightarrow 0,
$$
and
$$ 0 \longrightarrow F \longrightarrow (-,
X[1]) \stackrel{(-, u[1])}{\longrightarrow} (-, Y[1])
\stackrel{(-, v[1])}{\longrightarrow} (-, Z[1]) \stackrel{(-,
w[1])}{\longrightarrow} (-, X[2]) \longrightarrow \cdots $$

\vskip 10pt

Denote by ${\rm Coho}(\mathcal{T})$ the full subcategory of
$(\mathcal{T}^{\rm op}, {\rm Ab})$ consisting of cohomological
functors. By Lemmas \ref{lem:coherent} and
\ref{lem:trianglecoherent}, we deduce that the category
$\widehat{\mathcal{T}}$ is an abelian category which has enough
projective objects. In particular, we may define the extension
groups ${\rm Ext}^i(F, G)$ for any coherent functors $F$ and $G$,
$i\geq 1$.\par \vskip 5pt

The next lemma is also known; compare \cite[p.258]{Ne}.

\begin{lem}\label{lem:coho}
Let $H\in \widehat{\mathcal{T}}$. Then we have $H\in {\rm
Coho}(\mathcal{T})$ if and only if ${\rm Ext}^i(F, H)=0$ for all
$i\geq 1$ and $F\in \widehat{\mathcal{T}}$.
\end{lem}

\begin{proof} Assume that $F \in \widehat{\mathcal{T}}$. Note that we just
obtained a projective resolution for $F$. We take the resolution to
compute ${\rm Ext}^i(F, H)$. Apply Yoneda Lemma $((-, C), H)\simeq
H(C)$. We obtain that ${\rm Ext}^i(F, H)$ is just the $i$-th
cohomology group of the following complex
$$0 \longrightarrow H(Z) \longrightarrow H(Y) \longrightarrow H(X)
 \longrightarrow H(Z[-1]) \longrightarrow \cdots $$
Here we use the notation as in the above two-sided resolution. 
Now if $H$ is cohomological, the above complex is exact at all the
positions other than the zeroth one (where $H(Z)$ sits). Hence ${\rm
Ext}^i(F, H)=0$ for  all $i\geq 1$. \par

 On the other hand, assume that
${\rm Ext}^i(F, H)=0$ for all $i\geq 1$ and $F\in
\widehat{\mathcal{T}}$. To see that $H$ is cohomological, take a
distinguished triangle $X\stackrel{u}{\rightarrow} Y
\stackrel{v}{\rightarrow} Z \stackrel{w}{\rightarrow} X[1].$ Take
$F={\rm Cok}(-, v)$; it is a coherent functor. As in the above
discussion, the assumption that ${\rm Ext}^i(F, H)=0$ will imply
that the complex
$$0 \longrightarrow H(Z) \longrightarrow H(Y) \longrightarrow H(X)
 \longrightarrow H(Z[-1]) \longrightarrow \cdots $$
is exact at all the positions other than the zeroth one. This
implies that $H$ is cohomological. \end{proof}

Denote by ${\rm add}({\bf p}(\mathcal{T}))$ the full subcategory of
$(\mathcal{T}^{\rm op}, {\rm Ab})$ consisting of direct summands of
representable functors. Observe that ${\rm add}({\bf
p}(\mathcal{T}))\subseteq {\rm Coho}(\mathcal{T})\cap
\widehat{\mathcal{T}}$.

\vskip 5pt

 The following result is due to Freyd and Verdier
independently (\cite{Fr, V}). We refer to \cite[p.169]{Ne} and
\cite[4.2]{Kr} for a modern treatment. The result we present is in a
slightly different form. Recall that an abelian category is
\emph{Frobenius} provided that it has enough projective and enough
injective objects, and the class of projective objects coincides
with the class of injective objects; see \cite[Chapter I.2]{H1}.

\begin{prop}\label{prop:Freyd}
 Let $\mathcal{T}$ be a triangulated category. Then the category $\widehat{\mathcal{T}}$ of coherent
 functors is a Frobenius abelian category  such that the class of its projective-injective objects
equals
$${\rm Coho}(\mathcal{T}) \cap
\widehat{\mathcal{T}}={\rm add}({\bf p}(\mathcal{T})).$$
\end{prop}

\begin{proof} It is already observed that $\widehat{\mathcal{T}}\subseteq
(\mathcal{T}^{\rm op}, {\rm Ab})$ is an abelian subcategory which
has enough projective objects. The class of projective objects
equals ${\rm add}({\bf p}(\mathcal{T}))$. Note that $(-, C)$ is a
cohomological functor. By Lemma \ref{lem:coho}, we have, for each $F
\in \widehat{\mathcal{T}}$ and $C \in \mathcal{T}$, ${\rm Ext}^i(F,
(-, C))=0$, $i\geq 1$. Therefore, the representable functor $(-, C)$
is an injective object in $\widehat{\mathcal{T}}$. Hence we have
shown that $\widehat{\mathcal{T}}$ has enough projective objects and
that every projective object is injective.
\par

We claim that $\widehat{\mathcal{T}}$ has enough injective objects
and that every injective object is projective. To see this, assume
that $F$ is a coherent functor, and take a presentation $(-,
Y)\stackrel{(-, v)}{\rightarrow} (-, Z) \rightarrow F \rightarrow
0$. Then we have seen that there is a monomorphism $F \rightarrow
(-, X[1])$ via the two-sided resolution of $F$. By the above proof,
the functor $(-,X[1])$ is injective. Then the category
$\widehat{\mathcal{T}}$ has enough injective objects. Assume further
that $F$ is an injective object in $\widehat{\mathcal{T}}$. Then the
monomorphism $F \rightarrow (-, X[1])$ splits, and hence $F$ is
projective. This proves the claim, and then we have shown that
$\widehat{\mathcal{T}}$ is a Frobenius category and the class of its
projective-injective equals ${\rm add}({\bf p}(\mathcal{T}))$.  \par

It remains to show ${\rm Coho}(\mathcal{T}) \cap
\widehat{\mathcal{T}}={\rm add}({\bf p}(\mathcal{T}))$. We have
already noticed that ${\rm add}({\bf p}(\mathcal{T}))\subseteq {\rm
Coho}(\mathcal{T}) \cap \widehat{\mathcal{T}}$. Let $H \in {\rm
Coho}(\mathcal{T}) \cap \widehat{\mathcal{T}} $. By Lemma
\ref{lem:coho}, ${\rm Ext}^i(F, H)=0$ for all $F \in
\widehat{\mathcal{T}}$ and $i\geq 0$. Then $H$ is an injective
object in $\widehat{\mathcal{T}}$, and hence by above $H$ lies in
${\rm add}({\bf p}(\mathcal{T}))$.
\end{proof}

 In what follows $R$ will be a commutative artinian ring. We denote by $R\mbox{-mod}$ the
 category of finitely generated $R$-modules.  The
triangulated category $\mathcal{T}$ will be $R$-linear which is
Hom-finite  and Krull-Schmidt. Consider the category
$(\mathcal{T}^{\rm op}, R\mbox{-mod})$ of $R$-linear functors from
$\mathcal{T}$ to $R\mbox{-mod}$. Observe that $(\mathcal{T}^{\rm
op}, R\mbox{-mod})$ is viewed as a full subcategory of
$(\mathcal{T}^{\rm op},{\rm  Ab})$ via the forgetful functor. Recall
the notion of generalized Serre duality from the introduction.

\begin{prop}\label{prop:domainrange}
The full subcategories $\mathcal{T}_r$ and $\mathcal{T}_l$ are thick
triangulated subcategories of $\mathcal{T}$.
\end{prop}

\begin{proof} We only prove the result on $\mathcal{T}_r$. We first
claim that $$\mathcal{T}_r=\{X\in \mathcal{T}\;  | \; D(X, -) \in
\widehat{\mathcal{T}}\}.$$ To see this, note that $\mathcal{T}$ is
Krull-Schmidt, hence it is idempotent-split; see \cite{H1} and
\cite[Theorem A.1]{CYZ}. Then we have ${\rm add}({\bf
p}(\mathcal{T}))={\bf p}(\mathcal{T})$. Note that the functor
$D(X,-)$ is always cohomological. Hence if $D(X,-)$ is coherent,
then by Proposition \ref{prop:Freyd}, $D(X, -)\in {\rm
Coho}(\mathcal{T})\cap \widehat{\mathcal{T}}={\rm add}({\bf
p}(\mathcal{T}))={\bf p}(\mathcal{T})$, that is, $D(X, -)$ is
representable. Hence, the functor $D(X, -)$ is coherent if and only
if it is representable. Now the claim follows.\par

Note that $D(X[n], -) \simeq D(X, -)\circ [-n]$; here for each
integer $n$, $[n]$ denotes the $n$-th power of the translation
functor $[1]$ on $\mathcal{T}$. Then the subcategory
$\mathcal{T}_r\subseteq \mathcal{T}$ is closed under the functors
$[n]$. Let $X'\rightarrow X \rightarrow X''\rightarrow X'[1]$ be a
distinguished triangle such that  $X', X'' \in \mathcal{T}_r$. Then
we have an exact sequence of functors
\begin{align*}
(X'[1], -)\longrightarrow (X'', -) \longrightarrow (X, -)
\longrightarrow (X', -)\longrightarrow (X''[-1], -),
\end{align*}
and hence the following sequence
\begin{align*}
D(X''[-1], -)\longrightarrow D(X', -) \longrightarrow D(X, -)
\longrightarrow D(X'', -)\longrightarrow D(X'[1], -)
\end{align*}
is also exact. Since $X', X'' \in \mathcal{T}_r$, the four terms
other than $D(X, -)$ in the above sequence lie in
$\widehat{\mathcal{T}}$. By Lemmas \ref{lem:coherent} and
\ref{lem:trianglecoherent} $\widehat{\mathcal{T}}\subseteq
(\mathcal{T}^{\rm op}, {\rm  Ab})$ is an abelian subcategory which
is closed under extensions. Hence we infer that $D(X, -) \in
\widehat{\mathcal{T}}$, and then by the claim above, we have $X \in
\mathcal{T}_r$. We have shown that $\mathcal{T}_r$ is a triangulated
subcategory of $\mathcal{T}$. Observe that $\mathcal{T}_r$ is thick,
that is, it is closed under taking direct summands, since direct
summands of a coherent functor are still coherent.
\end{proof}

\par \vskip 10pt

Let us assume  that $\mathcal{C}$ is  a Hom-finite $R$-linear
additive category which  is \emph{skeletally small}, that is, the
iso-classes of objects in $\mathcal{C}$ form a set.  Consider the
full subcategory $(\mathcal{C}^{\rm op}, R\mbox{-mod})$ of
$(\mathcal{C}^{\rm op},{\rm  Ab})$ consisting of $R$-linear
functors. Then we have duality functors \begin{center}$D\colon
(\mathcal{C}^{\rm op}, R\mbox{-mod}) \longrightarrow (\mathcal{C},
R\mbox{-mod})$ and $D\colon (\mathcal{C}^{\rm op}, R\mbox{-mod})
\longrightarrow (\mathcal{C}, R\mbox{-mod})$,\end{center} which are
induced by the Matlis duality on $R\mbox{-mod}$. The additive
category $\mathcal{C}$ is called a {\em dualizing $R$-variety} if
the two duality functors $D$ preserve coherent functors (\cite{AR}).
It is easy to see that $\mathcal{C}$ is a dualizing $R$-variety if
and only if $\mathcal{C}$ has pseudokernels and pseudocokernels, and
$D(X, -) \in \widehat{\mathcal{C}}$ and $D(-, X) \in
\widehat{\mathcal{C}^{\rm op}}$ for each object $X\in \mathcal{C}$; 
one may apply Lemma \ref{lem:coherent} and  \cite[Theorem 2.4(1)]{AR}.

Recall from Appendix A.1 that a Hom-finite category $\mathcal{C}$ is
said to have Serre duality, if
$\mathcal{C}_r=\mathcal{C}=\mathcal{C}_l$. Combining Lemma
\ref{lem:trianglecoherent} with the claim (and its dual) in the
above proof, we get the following observation; compare
\cite[Proposition 2.11]{IY}.

\begin{cor}
Let $\mathcal{T}$ be a skeletally small Hom-finite Krull-Schmidt
triangulated category. Then $\mathcal{T}$ has Serre duality if and
only if $\mathcal{T}$ is a dualizing $R$-variety. \hfill $\square$
\end{cor}

The following result generalizes slightly a result due to Bondal and
Kapranov (\cite[Proposition 3.3]{BK}). Let us stress that the tricky
proof is modified from an argument independently due to Huybrechts
(\cite[Proposition 1.46]{Huy}) and Van den Bergh (\cite[Theorem
A.4.4]{Bo}). Recall the notions of generalized Serre duality and
trace function in the appendix. For the notion of triangle functor,
we refer to \cite[section 8]{Ke}.

\begin{prop}\label{prop:2.7}
Let $\mathcal{T}$ a Hom-finite Krull-Schmidt triangulated category.
Then there is a natural isomorphism $\eta_X\colon S(X[1])\rightarrow
S(X)[1]$ for each $X \in \mathcal{T}_r$ such that the pair $(S,
\eta)$ is a triangle  functor between the triangulated categories
$\mathcal{T}_r$ and $\mathcal{T}_l$.
\end{prop}

\begin{proof} Recall from Proposition \ref{prop:domainrange} that both
$\mathcal{T}_r$ and $\mathcal{T}_l$ are triangulated categories of
$\mathcal{T}$. We will construct a natural  isomorphism
$$\eta_X\colon S(X[1]) \longrightarrow S(X)[1] $$
for each $X\in \mathcal{T}_r$, such that the pair $(S, \eta)$ is a
triangle  functor.\par

Let $X\in \mathcal{T}_r$ and $Y \in \mathcal{T}$. Consider the
following composite of natural $R$-linear isomorphisms
\begin{align*}
\Psi_{X,Y}\; \colon \; (Y, S(X[1])) \stackrel{\phi_{X[1],
Y}^{-1}}{\longrightarrow}
D(X[1], Y) \stackrel{-D\circ [-1]}{\longrightarrow} D(X, Y[-1])\\
\stackrel{\phi_{Y[-1], X}}{\longrightarrow} (Y[-1], S(X))
\stackrel{[1]}{\longrightarrow} (Y, S(X)[1]),
\end{align*}
where $\phi$ is the natural isomorphism in Proposition
\ref{prop:A.1}. Note the minus sign in the second isomorphism above.
By following the composition carefully, we obtain
that, for each $f\in (X, Y[-1])$ and $g\in (Y, S(X[1]))$,
$$(f,\;  \Psi_{X,Y}(g)[-1])_{X,Y[-1]}=-(f[1], \; g)_{X[1], Y}.$$

Note that $\Psi$ defines an isomorphism $\Psi_X\colon (-, S(X[1]))
\simeq (-, S(X)[1])$ of $R$-linear functors. By Yoneda Lemma there
is a unique isomorphism $\eta_X\colon S(X[1])\rightarrow S(X)[1]$
such that $\Psi_X= (-, \eta_X)$, in other words,
$\Psi_{X,Y}(g)=\eta_X\circ g$. Consequently, we have
\begin{align}
(f, (\eta_X\circ g)[-1])_{X,Y[-1]}=-(f[1], g)_{X[1], Y}.
\end{align}
We observe that $\Psi_{X, Y}$ is natural both in $X$ and $Y$. It
follows then that $\eta_X$ is natural in $X$. In other words, the
isomorphisms $\eta_X$ define a natural isomorphism of functors.

\vskip 3pt

 Next we show that the pair $(S, \eta)$ is a triangle functor.
Given a distinguished triangle $X\stackrel{u}{\rightarrow} Y
\stackrel{v}{\rightarrow} Z \stackrel{w}{\rightarrow} X[1]$ in
$\mathcal{T}_r$. We have to show that the following triangle is also
distinguished (in $\mathcal{T}_l$, or in $\mathcal{T}$)
$$S(X)\stackrel{S(u)}{\longrightarrow} S(Y)
\stackrel{S(v)}{\longrightarrow} S(Z) \stackrel{\eta_X\circ
S(w)}{\longrightarrow} S(X)[1].$$ \vskip 3pt

 We claim that
the following sequence of functors  is exact
$$(-, S(Y)) \longrightarrow (-, S(Z))\longrightarrow (-, S(X)[1])
 \longrightarrow (-, S(Y)[1]) \longrightarrow (-, S(Z)[1]).$$
In fact, from the distinguished triangle $X\stackrel{u}{\rightarrow}
Y \stackrel{v}{\rightarrow} Z \stackrel{w}{\rightarrow} X[1]$, we
infer that the following sequence is exact
$$(Z[1], -) \longrightarrow (Y[1], -) \longrightarrow (X[1], -) \longrightarrow (Z,-) \longrightarrow (Y, -),$$
and hence the following sequence is also exact
$$D(Y, -) \longrightarrow D(Z, -) \longrightarrow D(X[1], -) \longrightarrow D(Y[1],-) \longrightarrow D(Z[1], -).$$
Note that in the following commutative diagram of functors all
vertical morphisms are isomorphisms
\[\xymatrix{
D(Y, -) \ar[d]^-{\phi_Y} \ar[r] & D(Z, -) \ar[d]^-{\phi_Z} \ar[r]
&D(X[1],-) \ar[d]^-{\Psi_X\circ \phi_{X[1]}} \ar[r] &D(Y[1],
-)\ar[d]^-{\Psi_Y\circ \phi_{Y[1]}} \ar[r]
&D(Z[1],-) \ar[d]^-{\Psi_Z\circ \phi_{Z[1]}}\\
(-, S(Y)) \ar[r] &(-, S(Z)) \ar[r] &(-, S(X)[1]) \ar[r] &(-,
S(Y)[1]) \ar[r] & (-, S(Z)[1]).
 }\]
Here $\phi_{X[1]}$, $\phi_{Y[1]}$ and $\phi_{Z[1]}$ are the natural
isomorphisms of functors induced by $\phi$ in Proposition
\ref{prop:A.1}. This proves the claim.\par \vskip 3pt

Take a distinguished triangle $S(X)\stackrel{S(u)}{\rightarrow} S(Y)
\stackrel{\alpha}{\rightarrow} W \stackrel{\beta}{\rightarrow}
S(X)[1]$. We claim that to prove the result, it suffices to show
that there is a morphism $\delta\colon W \rightarrow S(Z)$ such that
the following diagram commutes.

\[\xymatrix{
S(X) \ar@{=}[d] \ar[r]^-{S(u)} & S(Y) \ar@{=}[d] \ar[r]^-{\alpha}
& W \ar[d]^-{\delta} \ar[r]^-\beta & S(X)[1] \ar@{=}[d]\\
S(X) \ar[r]^-{S(u)} & S(Y) \ar[r]^-{S(v)}  & S(Z)
 \ar[r]^-{\eta_X\circ S(w)} & S(X)[1]}\]
In fact, assume that we already have the above commutative diagram.
Then we have the following commutative diagram of functors
\[\xymatrix{
(-, S(X)) \ar[r] \ar@{=}[d] & (-, S(Y)) \ar[r] \ar@{=}[d] & (-, W)
\ar[r] \ar[d]^{(-, \delta)} & (-, S(X)[1]) \ar[r] \ar@{=}[d] &(-,
S(Y)[1]) \ar@{=}[d]\\
(-, S(X)) \ar[r] &(-, S(Y)) \ar[r] &(-, S(Z)) \ar[r] &(-, S(X)[1])
\ar[r] &(-, S(Y)[1]) }\]
Note the the rows are exact (see the
previous claim). Then Five Lemma implies that $(-, \delta)$ is an
isomorphism, and hence by Yoneda Lemma $\delta$ is an isomorphism.
We are done.

\par\vskip 3pt

To complete the proof, it suffices to find a required morphism
$\delta$. To find such a morphism $\delta\in (W, S(Z))$ is to solve
the two equations $\delta\circ\alpha=S(v)$ and $\eta_X\circ
S(w)\circ \delta=\beta$. By (2.1) we have
\begin{align}
{\rm Tr}_X ((\eta_X \circ f)[-1]) = - {\rm Tr}_{X[1]} (f), \quad
f\in (X[1], S(X[1])).
\end{align}
For the definition of Tr and the  bilinear form $(-, -)$, see
Appendix. By the non-degeneratedness of the bilinear form $(-, -)$,
we infer that $\delta\circ\alpha=S(v)$ is equivalent to the
equations ${\rm Tr_Z}(S(v)\circ x)={\rm Tr}_Z(\delta\circ \alpha
\circ x)$ for all $x\in (Z, S(Y))$; and that $\eta_X\circ S(w)\circ
\delta=\beta$ is equivalent to the equations ${\rm
Tr}_X((\eta_X\circ S(w)\circ \delta)[-1]\circ y)={\rm
Tr}_X(\beta[-1]\circ y)$ for all $y\in (X, W[-1])$. \par

Note that we have
\begin{align*}
{\rm Tr}_X((\eta_X\circ S(w)\circ \delta)[-1]\circ y) &= - {\rm
Tr}_{X[1]} (S(w)\circ \delta \circ y[1]) \\
 &=-{\rm Tr}_Z (\delta
\circ y[1] \circ w),
\end{align*}
where the first equality uses (2.2) and the second uses (A.4). Hence
it suffices to find $\delta \in (W, S(Z))$ such that
\begin{align*}
{\rm Tr}_Z(\delta\circ \alpha\circ x)&= {\rm Tr}_Z(S(v)\circ x),\quad \forall\;  x \in (Z, S(Y)),\\
{\rm Tr}_Z(\delta \circ y[1] \circ w)&=-{\rm Tr}_X(\beta[-1]\circ
y), \quad \forall\; y\in (X, W[-1]).
\end{align*}

By  Proposition \ref{prop:A.1} and its proof  we have the
isomorphism $\phi_{W, Z}\colon (W, S(Z)) \simeq D(Z, W)$ and the
equality $\phi_{W, Z}(\delta)={\rm Tr}_Z(\delta\circ -)$. Here we
recall that $E$ is the minimal injective cogenerator for $R$ and
that $D={\rm Hom}_R(-, E)$. Hence to complete the proof, it suffices
to find an $R$-linear morphism $F$ from $(Z, W)$ to $E$ such that
\begin{align*}
F( \alpha\circ x)&= {\rm Tr}_Z(S(v)\circ x),\quad \forall\; x \in (Z, S(Y)),\\
F( y[1] \circ w)&=-{\rm Tr}_X(\beta[-1]\circ y), \quad \forall\; y
\in (X, W[-1]).
\end{align*}
Using the injectivity of $E$, it is not hard to see that such a
morphism $F$ exists provided that whenever $\alpha\circ x=y[1]\circ
w$, then ${\rm Tr}_Z(S(v)\circ x)=-{\rm Tr}_X(\beta[-1]\circ y)$.
Now we assume that $\alpha\circ x=y[1]\circ w$. Then  we have the
following morphism of distinguished triangles.

\[\xymatrix{
Y \ar[r]^-v \ar@{-->}[d]^{\phi_0} & Z \ar[d]^-x \ar[r]^-w &X[1] \ar[d]^-{y[1]} \ar[r]^-{-u[1]} & Y[1] \ar@{-->}[d]^{\phi_0[1]}\\
S(X) \ar[r]^-{S(u)} & S(Y) \ar[r]^-\alpha & W \ar[r]^-\beta &
S(X)[1] }\]

Therefore we have
\begin{align*}
{\rm Tr}_Z(S(v) \circ x)&= {\rm Tr}_Y(x \circ v)\\
                        &= {\rm Tr}_Y( S(u)\circ \phi_0)\\
                        &={\rm Tr}_X(\phi_0\circ u)\\
                        &={\rm Tr}_X(-\beta[-1]\circ y)\\
                        &=-{\rm Tr}_X(\beta[-1]\circ y),
\end{align*}
where the first and third equality use (A.4) and the second and
fourth use the commutativity of the diagram above. This completes
the proof.\end{proof}

Let $\mathcal{T}$ be a Hom-finite Krull-Schmidt triangulated
category as above. Recall from \cite{H1} and \cite[1.2]{H2} that an
\emph{Auslander-Reiten triangle} in $\mathcal{T}$ is a distinguished
triangle $X \stackrel{u}{\rightarrow} Y \stackrel{v}{\rightarrow} Z
\stackrel{w}{\rightarrow} X[1]$ such that both $X$ and $Z$ are
indecomposable, and that $w\neq 0$ satisfies that for each
non-retraction $\gamma\colon Y'\rightarrow Z$ in $\mathcal{T}$ we
have $w\circ \gamma$=0. By \cite[p.31]{H1} the morphism $u$ is left
almost-split, that is, any non-section $h\colon X \rightarrow Y'$
factors through $u$; dually, the morphism $v$ is right almost-split.

\par \vskip 5pt

 The following result is a slight generalization of a result  due to Reiten and Van den Bergh
 (\cite[Proposition I.2.3]{RV}), where the result is formulated in a different terminology.

\begin{prop} \label{prop:2.8}
Let $\mathcal{T}$ be a Hom-finite Krull-Schmidt triangulated category, and let $X \in \mathcal{T}$ be
an indecomposable object. Then we have
\begin{enumerate}
\item $X \in\mathcal{T}_r$ if and only if there is an Auslander-Reiten
triangle with  $X$ as the third term;
\item $X \in \mathcal{T}_l$ if
and only if there is an Auslander-Reiten triangle with $X$ as the
first term.\end{enumerate}
\end{prop}

\begin{proof} We only prove the first statement. By a dual argument, one
proves the second one.\par

For the ``only if" part, assume that $X \in \mathcal{T}_r$ is
indecomposable. Recall from Proposition \ref{prop:A.1} that the
generalized Serre functor $S\colon \mathcal{T}_r \rightarrow
\mathcal{T}_l$ is an equivalence. Therefore $S(X)$ is also
indecomposable.  For each $X\in \mathcal{T}_r$ and $Y \in
\mathcal{T}$ there is a non-degenerated bilinear form $(-, -)_{X,
Y}: (X, Y) \times (Y, S(X)) \rightarrow E$; see the proof of
Proposition \ref{prop:A.1}. Take a nonzero morphism $w\colon X
\rightarrow S(X)$ such that $({\rm rad}\; {\rm End}_\mathcal{T}(X),
w)_{X, X}=0$. Here ${\rm rad}\;{\rm End}_\mathcal{T}(X)$ denotes the
Jacobson radical of ${\rm End}_\mathcal{T}(X)$. Then we form a
distinguished triangle $S(X)[-1]\stackrel{u}{\rightarrow} Y
\stackrel{v}{\rightarrow} X \stackrel{w}{\rightarrow} S(X)$. We
claim that it is an Auslander-Reiten triangle. Then we are done.

In fact, it suffices to show that each morphism $\gamma\colon
X'\rightarrow X$, which is not split-epic, satisfies  $w\circ
\gamma=0$. To see this, for any $x \in (X, X')$, consider $(x,
w\circ \gamma)_{X, X'}=(\gamma \circ x, w)_{X, X}$; see (A.2).
Since $\gamma$ is a not split-epic, then $\gamma\circ x \in {\rm rad}\;{\rm
End}_\mathcal{T}(X)$, and thus $(\gamma\circ x, w)_{X, X}=0$.
Therefore $(x, w\circ \gamma)_{X, X'}=0$. By the non-degeneratedness
of the bilinear form, we have $w \circ \gamma=0$.\par\vskip 3pt

For the ``if" part, assume that we have an Auslander-Reiten triangle
$Z\stackrel{u}{\rightarrow} Y \stackrel{v}{\rightarrow} X
\stackrel{w}{\rightarrow} Z[1]$. Take any $R$-linear morphism ${\rm
Tr}\colon (X, Z[1]) \rightarrow E$ such that ${\rm Tr}(w)\neq 0$.
Here we use that $E$ is an injective cogenerator. We claim that for
each $X'\in \mathcal{T}$ the pairing
$$(-, -):(X, X') \times (X' , Z[1])
\longrightarrow E$$ given by $(f, \; g)={\rm Tr}(g\circ f)$ is
non-degenerated. If so, we have an induced isomorphism of $R$-linear
functors $D(X, -) \simeq (-, Z[1])$. Hence $X\in \mathcal{T}_r$, and
then we are done.
\par

To prove the claim, let $X' \in \mathcal{T}$. Let $f\colon X
\rightarrow X'$ be a nonzero morphism, and let $X''
\stackrel{s}{\rightarrow} X \stackrel{f}{\rightarrow} X' \rightarrow
X''[1]$ be a distinguished triangle. Since $f \neq 0$, then $s$ is
not split-epic. Then by the properties of Auslander-Reiten triangle,
we have $w\circ s=0$. This implies that $w$ factors through $f$,
that is, there is a morphism $g\colon X'\rightarrow Z[1]$ such that
$w=g\circ f$, and hence $(f,\; g)\neq 0$. On the other hand, let
$g\colon X' \rightarrow Z[1]$ be a nonzero morphism, and let $Z
\stackrel{u'}{\rightarrow} X''' \rightarrow X'
\stackrel{g}{\rightarrow} Z[1]$ be a distinguished triangle. Since $g\neq 0$, then $u'$ is not
split-mono. Then by the properties of Auslander-Reiten triangle, we
deduce that $u'$ factors thorough $u$, and then we have the
following morphism of distinguished triangles

\[\xymatrix{
Z \ar@{=}[d] \ar[r]^-u & Y \ar@{-->}[d] \ar[r]^-v & X
\ar@{-->}[d]^-f \ar[r]^-w & Z[1] \ar@{=}[d]\\
Z \ar[r]^{u'} & X''' \ar[r] &X' \ar[r]^-g & Z[1] }\]
 Hence $g\circ f=w$ and then $(f,\;
g)\neq 0$. This completes the proof. \end{proof}

\section{Computing generalized Serre duality}

In this section, we will compute explicitly the generalized Serre
duality for the bounded derived categories of artin algebras and of
certain noncommutative projective schemes in the sense of Artin and
Zhang (\cite{AZ}); see Theorem \ref{thm:3.5} and Theorem
\ref{thm:3.10}.

\vskip 5pt

Let $R$ be a commutative artinian ring. Let $\mathcal{C}$ be an
$R$-linear category which is Hom-finite and which is assumed to be
additive. Denote by $S_\mathcal{C}\colon \mathcal{C}_r \rightarrow
\mathcal{C}_l$ its generalized Serre functor. Denote by $(-, -)_{X,
Y}\colon (X, Y)\times (Y, S_{\mathcal{C}}(X))\rightarrow E$ the
associated bilinear form for each $X \in \mathcal{C}_r$ and $Y\in
\mathcal{C}$. Here $E$ is the minimal injective cogenerator of $R$.
For details, see Appendix.

Let $\mathbf{K}^b(\mathcal{C})$ be the bounded homotopy category of
complexes in $\mathcal{C}$. Note that the category
$\mathbf{K}^b(\mathcal{C})$ is naturally $R$-linear and then is
Hom-finite. Denote its generalized Serre functor by $S\colon
\mathbf{K}^b(\mathcal{C})_r\rightarrow \mathbf{K}^b(\mathcal{C})_l$.
The following result is inspired by the argument in \cite[p.37]{H1}.

\begin{prop}\label{prop:3.1}
Use the notation as above. Then we have
 \begin{center} $\mathbf{K}^b(\mathcal{C}_r)
\subseteq \mathbf{K}^b(\mathcal{C})_r$, $\mathbf{K}^b(\mathcal{C}_l)
\subseteq \mathbf{K}^b(\mathcal{C})_l$ and
$S|_{\mathbf{K}^b(\mathcal{C}_r)}=S_\mathcal{C}$,\end{center} where
$S_\mathcal{C}$ is viewed as a functor on complexes between
$\mathbf{K}^b(\mathcal{C}_r)$ and $\mathbf{K}^b(\mathcal{C}_l)$.
\end{prop}

\begin{proof} To prove the result, it suffices to build up a natural
non-degenerated bilinear form
\begin{align*}
(-, -)_{X^\bullet, Y^\bullet} \colon {\rm
Hom}_{\mathbf{K}^b(\mathcal{C})}(X^\bullet, Y^\bullet)\times {\rm
Hom}_{\mathbf{K}^b(\mathcal{C})}(Y^\bullet,
S_\mathcal{C}(X^\bullet))\longrightarrow E
\end{align*}
for each $X^\bullet \in \mathbf{K}^b(\mathcal{C}_r)$ and $Y^\bullet
\in \mathbf{K}^b(\mathcal{C})$. We define
$$(f^\bullet, g^\bullet)_{X^\bullet,
Y^\bullet}=\sum_{i\in \mathbb{Z}}(-1)^i(f^i, g^i)_{X^i, Y^i},$$
where $f^\bullet=(f^i)$ and $g^\bullet=(g^i)$ are chain
morphisms.\par

Note that the bilinear form is well defined. In fact, if $f^\bullet$
is homotopic to $0$, then $(f^\bullet, g^\bullet)_{X^\bullet,
Y^\bullet}=0$. Assume that $f^i=d_Y^{i-1}\circ h^i+ h^{i+1}\circ
d_X^i$ for each $i\in \mathbb{Z}$. Here $\{h^i\colon X^i\rightarrow
Y^{i-1}\}_{i\in \mathbb{Z}}$ is the homotopy. Then we have
\begin{align*}
(f^\bullet, \; g^\bullet)_{X^\bullet, Y^\bullet}&=\sum_{i\in
\mathbb{Z}} (-1)^i(d_Y^{i-1}\circ h^i,\; g^i)_{X^i, Y^i} + (-1)^i
(h^{i+1}\circ
d_X^i,\; g^i)_{X^i, Y^i}\\
&= \sum_{i\in \mathbb{Z}} (-1)^i (h^i,\; g^i \circ d_Y^{i-1})_{X^i,
Y^{i-1}} + (-1)^i  (h^{i+1}, \; S_\mathcal{C}(d_X^i)\circ
g^i)_{X^{i+1}, Y^i}\\
&=\sum_{i \in \mathbb{Z}} (-1)^i (h^i, \; g^i\circ
d_Y^{i-1}-d_{S_{\mathcal{C}}(X)}^{i-1}\circ g^{i-1})_{X^i,
Y^{i-1}}\\
&=0.
\end{align*}
The second equality follows from the identities (A.1) and (A.2) in
the appendix, and the last one  follows from the fact that
$g^\bullet \colon Y^\bullet \rightarrow S_\mathcal{C}(X^\bullet)$ is
a chain morphism.

 To prove the naturalness of the bilinear form is  to verify the
following two identities
\begin{align*}
(f^\bullet \circ \theta^\bullet, \; g^\bullet)_{X'^\bullet,
Y^\bullet}&=(f^\bullet, \;S_\mathcal{C}(\theta^\bullet)\circ
g^\bullet)_{X^\bullet,
Y^\bullet}\\
(\gamma^\bullet\circ f^\bullet, \; g^\bullet)_{X^\bullet,
Y^\bullet}&= (f^\bullet, \; g^\bullet\circ
\gamma^\bullet)_{X^\bullet, Y'^\bullet}.
\end{align*}
Note that they can be easily derived from the identities (A.1) and
(A.2)  for the category $\mathcal{C}$.\par

We will prove the non-degeneratedness. For this end,  fix a complex
$X^\bullet \in \mathbf{K}^b(\mathcal{C}_r)$. Denote by
 $$\phi:
D{\rm Hom}_{\mathbf{K}^b(\mathcal{C})}(X^\bullet, -)\longrightarrow
{\rm Hom}_{\mathbf{K}^b(\mathcal{C})}(-, S_\mathcal{C}(X^\bullet))$$
the natural transformation induced from the above bilinear form. We
claim that $\phi$ is a natural isomorphism. Then we are done.

In fact, since $\phi$ is a natural transformation between two
cohomological functors, so by devissage, it suffices to prove that
$$\phi_{Y^\bullet}: D{\rm Hom}_{\mathbf{K}^b(\mathcal{C})}(X^\bullet,
Y^\bullet)\longrightarrow {\rm
Hom}_{\mathbf{K}^b(\mathcal{C})}(Y^\bullet,
S_\mathcal{C}(X^\bullet))$$
 is an isomorphism for each stalk complex
$Y^\bullet$. Without loss of generality, assume that $Y^\bullet=Y$
for some object $Y \in\mathcal{C}$. Then $D{\rm
Hom}_{\mathbf{K}^b(\mathcal{C})}(X^\bullet, Y)$ is just the zeroth
cohomology of the following complex
\begin{align*}
\cdots \longrightarrow D(X^{-1}, Y) \longrightarrow D(X^0, Y)
\longrightarrow D(X^1, Y) \longrightarrow \cdots
\end{align*}
By the generalized Serre duality of $\mathcal{C}$, the above complex
isomorphic to the following complex
\begin{align*}
\cdots \longrightarrow (Y, S_\mathcal{C}(X^{-1})) \longrightarrow (
Y, S_\mathcal{C}(X^0)) \longrightarrow (Y, S_\mathcal{C}(X^1))
\longrightarrow \cdots
\end{align*}
While the zeroth cohomology  of this complex is nothing but ${\rm
Hom}_{\mathbf{K}^b(\mathcal{C})}(Y, S_\mathcal{C}(X^\bullet))$. This
proves that $\phi_Y$ is an isomorphism.
\end{proof}

Let $\mathcal{A}$ be an $R$-linear abelian category.  Denote by
$\mathbf{D}^b(\mathcal{A})$ its bounded derived category, which  is
naturally $R$-linear. We will assume that
$\mathbf{D}^b(\mathcal{A})$  is Hom-finite. In this case the abelian
category $\mathcal{A}$ is said to be {\em Ext-finite} (\cite{LC}).
Observe that an Ext-finite abelian category is Hom-finite. By
\cite[Corollary 2.10]{BS} the category $\mathbf{D}^b(\mathcal{A})$
is idempotent-split, hence by \cite[Theorem A.2]{CYZ} it is
Krull-Schmidt; also see \cite[Corollary B]{LC}. Denote by
$\mathcal{P}$ and $\mathcal{I}$ the full subcategory of
$\mathcal{A}$ consisting of projective objects and injective
objects, respectively. Let $S_\mathcal{A}\colon \mathcal{A}_r
\rightarrow \mathcal{A}_l$ be the generalized Serre functor on
$\mathcal{A}$.

 We claim that
$\mathcal{A}_r\subseteq \mathcal{P}$ and $\mathcal{A}_l\subseteq
\mathcal{I}$. Take an object $X \in \mathcal{A}_r$. By the
generalized Serre duality, we have an isomorphism of functors $(X,
-) \simeq D(-, S_\mathcal{A}(X))$. Hence the functor $(-, X)$ is
right exact, and then the object $X$ is projective. Similarly each
object in $\mathcal{A}_l$ is injective. \vskip 5pt

 Note that we have natural isomorphisms
\begin{align*}
 {\rm Hom}_{\mathbf{D}^b(\mathcal{A})}(P^\bullet, X^\bullet)\simeq {\rm
Hom}_{\mathbf{K}^b(\mathcal{A})}(P^\bullet, X^\bullet), \quad
 {\rm Hom}_{\mathbf{D}^b(\mathcal{A})}(X^\bullet, I^\bullet)\simeq {\rm
Hom}_{\mathbf{K}^b(\mathcal{A})}(X^\bullet, I^\bullet),
\end{align*}
for any $P^\bullet \in \mathbf{K}^b(\mathcal{P})$, $X^\bullet\in
\mathbf{K}^b(\mathcal{A})$ and $I^\bullet \in
\mathbf{K}^b(\mathcal{I})$. In particular,
$\mathbf{K}^b(\mathcal{P})$ and $\mathbf{K}^b(\mathcal{I})$ are
naturally viewed as triangulated subcategories of
$\mathbf{D}^b(\mathcal{A})$; compare \cite[1.4]{H3} and
\cite[Example 12.2]{Ke}.

\vskip 5pt

 We have  the following direct consequence of Proposition
\ref{prop:3.1}.

\begin{cor}\label{cor:3.2}
Let $\mathcal{A}$ be an Ext-finite abelian category. Then we have
\begin{center} $\mathbf{K}^b(\mathcal{A}_r) \subseteq \mathbf{D}^b(\mathcal{A})_r$,
$\mathbf{K}^b(\mathcal{A}_l)\subseteq \mathbf{D}^b(\mathcal{A})_l$
and $S|_{\mathbf{K}^b(\mathcal{A}_r)}=S_\mathcal{A}$, \end{center}
where $S$
denotes the generalized Serre functor on
$\mathbf{D}^b(\mathcal{A})$.\hfill $\square$
\end{cor}

 \vskip 3pt

Following \cite[Appendix]{NV} a complex $X^\bullet$ in
$\mathbf{D}^b(\mathcal{A})$ is said to {\em have finite projective
dimension}, if there exists some integer $N$ such that ${\rm
Hom}_{\mathbf{D}^b(\mathcal{A})}(X^\bullet, M[n])=0$ for any $M \in
\mathcal{A}$ and $n \geq  N$. Set $\mathbf{D}^b(\mathcal{A})_{\rm
fpd}$ to be the full subcategory of $\mathbf{D}^b(\mathcal{A})$
consisting of complexes of finite projective dimension; it is a
thick triangulated subcategory.
\par \vskip 5pt

The following is well known; compare the proof of \cite[Proposition
6.2]{Ri}.

\begin{lem}\label{lem:3.3}
We have $\mathbf{K}^b(\mathcal{P}) \subseteq
\mathbf{D}^b(\mathcal{A})_{\rm fpd}$. On the other hand, if the
category $\mathcal{A}$ has enough projective objects, then
$\mathbf{K}^b(\mathcal{P})=\mathbf{D}^b(\mathcal{A})_{\rm fpd}.$
\hfill $\square$
\end{lem}

Dually one defines a complex of \emph{finite injective dimension}
and introduces the triangulated subcategory
$\mathbf{D}^b(\mathcal{A})_{\rm fid}$ of
$\mathbf{D}^b(\mathcal{A})$. Then one has the dual version of Lemma
\ref{lem:3.3}; compare \cite[Chapter I, Proposition 7.6]{Har}.

\vskip 5pt

We have the following observation.

\begin{prop}\label{prop:3.4}
Let $\mathcal{A}$ be an Ext-finite abelian category. Then we have
$\mathbf{D}^b(\mathcal{A})_r\subseteq \mathbf{D}^b(\mathcal{A})_{\rm
fpd}$ and $\mathbf{D}^b(\mathcal{A})_l \subseteq
\mathbf{D}^b(\mathcal{A})_{\rm fid}$.
\end{prop}

\begin{proof}
Take a complex $X^\bullet \in \mathbf{D}^b(\mathcal{A})_r$. Then we
have an isomorphism of functors over $\mathbf{D}^b(\mathcal{A})$:
${\rm Hom}_{\mathbf{D}^b(\mathcal{A})}(X^\bullet, -) \simeq D{\rm
Hom}_{\mathbf{D}^b(\mathcal{A})}(-, S(X^\bullet))$, where $S$ is the
generalized Serre functor on $\mathbf{D}^b(\mathcal{A})$. Take a
positive integer $N$ such that $H^n(S(X^\bullet))=0$ for any $n\leq
-N$; here $H^n(-)$ denotes the $n$-th cohomology of a complex. Let
$M \in \mathcal{A}$ and let $n\geq N$. Then ${\rm
Hom}_{\mathbf{D}^b(\mathcal{A})}(X^\bullet, M[n])\simeq D{\rm
Hom}_{\mathbf{D}^b(\mathcal{A})}(M[n], S(X^\bullet))=0$. Hence
$X^\bullet$ has finite projective dimension. Then we get the first
inclusion, and the second one can be proved dually. \end{proof}

\vskip 10pt

Here is our first explicit result on generalized Serre duality.
\vskip 5pt

\begin{thm}\label{thm:3.5}
Let $\mathcal{A}$ be a Hom-finite abelian category with enough
projective and enough injective objects. Assume that
$\mathcal{A}_r=\mathcal{P}$ and $\mathcal{A}_l=\mathcal{I}$. Denote
the generalized Serre functor $S_\mathcal{A}$ by $\nu$ which is
called the Nakayama functor. Then we have
$\mathbf{D}^b(\mathcal{A})_r=\mathbf{K}^b(\mathcal{P})$,
$\mathbf{D}^b(\mathcal{A})_l=\mathbf{K}^b(\mathcal{I})$ and the
generalized Serre functor on $\mathbf{D}^b(\mathcal{A})$ is given by
$\nu$ which operates on complexes term by term.
\end{thm}

\begin{proof}
 Note that under the assumption, the abelian category $\mathcal{A}$ is Ext-finite.
Then the result follows immediately from Corollary \ref{cor:3.2},
Lemma \ref{lem:3.3} and Proposition \ref{prop:3.4}.
\end{proof}

We have the following motivating example of Theorem \ref{thm:3.5}.

\begin{exm}\label{exm:3.6}{\rm (Happel)}
 Let $A$ be an artin $R$-algebra. Denote by $A\mbox{-{\rm mod}}$
 the category of finitely generated left $A$-modules. We denote by $A\mbox{-{\rm proj}}$
 (resp. $A\mbox{-{\rm inj}}$) the full subcategory consisting of projective (resp. injective) modules.

 It is well known that the category $A\mbox{-{\rm mod}}$ satisfies the assumption of
Theorem \ref{thm:3.5}. In this case, the functor $\nu$ is the usual
Nakayama functor $D{\rm Hom}_A(-, A)$; see \cite[p.37]{H1}. We apply
Theorem \ref{thm:3.5} to conclude that $\mathbf{D}^b(A\mbox{-{\rm
mod}})_r=\mathbf{K}^b(A\mbox{-{\rm proj}})$,
$\mathbf{D}^b(A\mbox{-{\rm mod}})_l=\mathbf{K}^b(A\mbox{-{\rm
inj}})$ and the generalized Serre functor is given by $\nu$.  We
apply this explicit computational result. We infer from the
statement (2) in Main Theorem  that in an Auslander-Reiten triangle
$X^\bullet\rightarrow Y^\bullet \rightarrow Z^\bullet \rightarrow
X^\bullet[1]$ in $\mathbf{D}^b(A\mbox{-{\rm mod}})$ we have
$X^\bullet\in \mathbf{K}^b(A\mbox{-{\rm inj}})$ and $Z^\bullet \in
\mathbf{K}^b(A\mbox{-{\rm proj}})$. Moreover, we have
$X^\bullet\simeq \nu(Z^\bullet)$; see the proof of Proposition
\ref{prop:2.8}. This observation is essentially due to Happel {\rm
(\cite[Theorem 1.4]{H2})}.
\end{exm}

In what follows we will consider the conditions (C) and (C') which
are introduced in Appendix A.2. We have the following observation.

\begin{lem}\label{lem:3.7}
Let $\mathcal{A}$ be a Hom-finite abelian category. Then we have
\begin{enumerate}
\item if $\mathcal{A}$ has enough projective objects and
$\mathcal{A}_r=\mathcal{P}$, then $\mathbf{D}^b(\mathcal{A})$
satisfies the condition (C);
\item if $\mathcal{A}$ has enough injective objects
and $\mathcal{A}_l=\mathcal{I}$, then $\mathbf{D}^b(\mathcal{A})$
satisfies the condition (C').
\end{enumerate}
\end{lem}

\begin{proof} We only prove (1). Assume that $\mathcal{A}$
has enough projective objects and that $\mathcal{A}_r=\mathcal{P}$.
By Corollary \ref{cor:3.2}, Lemma \ref{lem:3.3} and Proposition
\ref{prop:3.4}, we have
$\mathbf{D}^b(\mathcal{A})_r=\mathbf{K}^b(\mathcal{P})$.

Consider complexes $X^\bullet, {X^\bullet}' \in
\mathbf{D}^b(\mathcal{A})_r$ and $Z^\bullet \in
\mathbf{D}^b(\mathcal{A})$. Since $\mathcal{A}$ has enough
projective objects, we may take a quasi-isomorphism $P^\bullet \rightarrow Z^\bullet$
such that $P^\bullet$ is a bounded-above complex
of projective objects with finitely many cohomologies; see
\cite[Chapter I, Lemma 4.6]{Har}. Let $n\gg 0$. Take
${Z^\bullet}'=\sigma^{\geq -n}P^\bullet$ to be the brutal
truncation, and take $s$ to be the following composite
\begin{align*}
{Z^\bullet}' \longrightarrow P^\bullet \longrightarrow
Z^\bullet,
\end{align*}
where ${Z^\bullet}'\rightarrow P^\bullet$ is the natural chain map.
Observe that ${Z^\bullet}' \in \mathbf{D}^b(\mathcal{A})_r$. It is
direct to check that the following two maps induced by $s$
\begin{align*}
{\rm Hom}_{\mathbf{D}^b(\mathcal{A})}(X^\bullet, {Z^\bullet} ')
\rightarrow {\rm Hom}_{\mathbf{D}^b(\mathcal{A})}(X^\bullet,
Z^\bullet), \;  {\rm Hom}_{\mathbf{D}^b(\mathcal{A})}({Z^\bullet} ',
{X^\bullet} ') \rightarrow {\rm
Hom}_{\mathbf{D}^b(\mathcal{A})}(Z^\bullet, {X^\bullet} ')
\end{align*}
are isomorphisms. We are done with the condition (C). \end{proof}

\vskip 10pt

We will call an abelian category $\mathcal{A}$ a  {\em Gorenstein
category} provided that $\mathbf{D}^b(\mathcal{A})_{\rm fpd}$
coincides with $\mathbf{D}^b(\mathcal{A})_{\rm fid}$; compare
\cite{H3}. By combining Theorem \ref{thm:3.5}, Lemma \ref{lem:3.7}
and Proposition \ref{prop:A.4} together, we have the following
result.

\begin{prop}
Let $\mathcal{A}$ be a Hom-finite abelian category with enough
projective and enough injective objects. Assume that
$\mathcal{A}_r=\mathcal{P}$ and $\mathcal{A}_l=\mathcal{I}$. Then
the following statements are equivalent: \begin{enumerate}
\item the category $\mathbf{K}^b(\mathcal{P})$ has Serre duality;
\item the abelian category $\mathcal{A}$ is Gorenstein;
\item the category $\mathbf{K}^b(\mathcal{I})$ has Serre duality.
\hfill $\square$
\end{enumerate}
\end{prop}

The following immediate consequence gives a  characterization of
Gorenstein algebras, which seems to be of independent interest.
Recall that an artin $R$-algebra $A$ is \emph{Gorenstein} provided
that the regular module $A$ has finite injective dimension on both
sides (\cite{H3}). We observe that an artin algebra $A$ is
Gorenstein if and only if the abelian category $A\mbox{-mod}$ of
finitely generated left $A$-modules is Gorenstein in the sense
above.

\begin{cor}\label{cor:3.9} {\rm  (compare \cite[Theorem 3.4]{H3})}
Let $A$ be an artin $R$-algebra. Then the
following statements are equivalent:
\begin{enumerate}
\item the category $\mathbf{K}^b(A\mbox{-{\rm proj}})$ has Serre duality;
\item the artin algebra $A$ is  Gorenstein;
\item the category $\mathbf{K}^b(A\mbox{-{\rm inj}})$ has Serre
duality.\hfill $\square$
\end{enumerate}
\end{cor}

 \vskip 15pt

We will next consider the generalized Serre duality on the bounded
derived category of coherent sheaves on certain noncommutative
projective schemes in the sense of Artin and Zhang (\cite{AZ}).
\par \vskip 5pt

We  will follow \cite[section 2 and Appendix]{NV} closely. Recall
some standard notation. Let $k$ be a field. Let $A=\oplus_{n\geq 0}
A_n$ be a positively graded $k$-algebra which is \emph{connected},
that is, $A_0\simeq k$. We will assume that $A$ is left noetherian.
Denote by $A\mbox{-Gr}$ (\emph{resp}. $A\mbox{-gr}$) the category of
graded left $A$-modules (\emph{resp}. finitely-generated graded left
$A$-modules) with morphisms preserving degrees. For each $d\in
\mathbb{Z}$, denote by $(d)$ be the degree-shift functor on
$A\mbox{-Gr}$ and $A\mbox{-gr}$. For $j \geq 0$, denote by ${\rm
Ext}_{A\mbox{-}{\rm Gr}} ^j(-, -)$ and $\underline{\rm Ext}^j(-, -)$
the Ext functors on the category $A\mbox{-Gr}$ and the associated
graded Ext functors, respectively. Then we have  $\underline{\rm
Ext}^j(M, N)=\oplus_{d\in \mathbb{Z}} {\rm Ext}^j_{A\mbox{-}{\rm
Gr}}(M, N(d))$ for graded modules $M, N$. We will view $k$ as a
graded $A$-module concentrated at degree zero.

For a graded $A$-module $M$ denote by $\tau(M)$ the sum of all its
finite dimensional graded submodules. This gives rise naturally to
the torsion functor $\tau\colon A\mbox{-Gr} \rightarrow
A\mbox{-Gr}$. A graded module $M$ is called {\em a torsion module}
provided that $\tau(M)=M$. Denote by $A\mbox{-Tor}$ (\emph{resp}.
$A\mbox{-tor}$) the full subcategories of $A\mbox{-Gr}$
(\emph{resp}. $A\mbox{-gr}$) consisting of torsion modules. Note
that both of them are Serre subcategories. One defines the quotient
abelian categories
\begin{align*}
{\rm Tails}(A)=A\mbox{-Gr}/{A\mbox{-Tor}}\mbox{  and   }{\rm
tails}(A)=A\mbox{-gr}/{A\mbox{-tor}},
\end{align*}
which are called the category of {\em quasi-coherent} and {\em
coherent sheaves} on the noncommutative projective scheme ${\rm
Proj}(A)$, respectively. Observe that both categories are naturally
$k$-linear.

\par \vskip 5pt

A connected left noetherian algebra $A$ is said to \emph{satisfy the
condition ``$\chi$"} provided that for each module $M \in
A\mbox{-gr}$, the Ext spaces ${\underline{{\rm Ext}}}^j(k, M)$ are
finite dimensional for all  $j \geq 1$; compare \cite[Definition
3.2]{AZ}. In this case,  the abelian category ${\rm tails}(A)$ is
Ext-finite; see \cite[Proposition 2.2.1]{NV}.

\par \vskip 10pt

Here is our second explicit result on generalized Serre duality. Let
us remark that it is somehow a restatement of a result due to de
Naeghel and Van den Bergh (\cite[Theorem A.4]{NV}); compare
\cite[Proposition 7.48]{Ro} for  the commutative case. For the
notion of balanced dualizing complex and other unexplained notions,
we refer to \cite{Y} and \cite{NV}.

\begin{thm}\label{thm:3.10} Let $A$ be a connected  (two-sided) noetherian
algebra over a field $k$. Assume that
\begin{enumerate} \item the
algebra $A$ satisfies ``$\chi$"  and the functor $\tau$ is of finite
cohomological dimension; \item the opposite algebra $A^{\rm op}$
satisfies ``$\chi$" and the functor $\tau^{\rm op}$ is of finite
cohomological dimension.
\end{enumerate}
 Then we have  $\mathbf{D}^b({\rm
tails}(A))_r=\mathbf{D}^b({\rm tails}(A))_{\rm fpd}$,
$\mathbf{D}^b({\rm tails}(A))_l=\mathbf{D}^b({\rm tails}(A))_{\rm
fid}$ and the generalized Serre functor is given $S= -\otimes^{\bf
L}\mathcal{R}^\bullet[-1]$, where $\mathcal{R}^\bullet$ is induced
by the balanced dualizing complex  $R^\bullet$ of $A$.
\end{thm}

\begin{proof} Note that the assumption implies the existence of a balanced dualizing complex
$R^\bullet$ in the sense of \cite{Y}. Here the functor $
-\otimes^{\bf L}\mathcal{R}^\bullet$ is induced by the functor $-
\otimes^{\bf L} R^\bullet$. For details, see \cite[Appendix]{NV}.

Note that the abelian category ${\rm tails}(A)$ is Ext-finite. By
Proposition \ref{prop:3.4}, we have
\begin{center}
$\mathbf{D}^b({\rm tails}(A))_r\subseteq \mathbf{D}^b({\rm
tails}(A))_{\rm fpd}$ and $\mathbf{D}^b({\rm tails}(A))_l\subseteq
\mathbf{D}^b({\rm tails}(A))_{\rm fid}$.
\end{center}
We apply \cite[Theorem A.4]{NV}.  We obtain $\mathbf{D}^b({\rm
tails}(A))_r\supseteq \mathbf{D}^b({\rm tails}(A))_{\rm fpd}$;
moreover, the restriction of the generalized Serre functor
$S|_{\mathbf{D}^b({\rm tails}(A))_{\rm fpd}}$ is isomorphic to
$-\otimes^{\bf L}\mathcal{R}^\bullet[-1]$, which is an equivalence
between $\mathbf{D}^b({\rm tails}(A))_{\rm fpd}$ and
$\mathbf{D}^b({\rm tails}(A))_{\rm fid}$; here we implicitly use
\cite[Proposition 2.3]{AZ} and \cite[Lemma 2.2]{NV}. This
equivalence implies that $ \mathbf{D}^b({\rm tails}(A))_{\rm
fid}\subseteq \mathbf{D}^b({\rm tails}(A))_l $, completing the
proof.
\end{proof}
 \vskip 10pt

We end this section with two remarks.

\begin{rem} (1) Let $\mathcal{T}$ be a Hom-finite Krull-Schmidt
triangulated category. We remark that the Verdier quotient
triangulated categories $\mathcal{T}/{\mathcal{T}_r}$ and
$\mathcal{T}/{\mathcal{T}_l}$ measure how far the triangulated
category $\mathcal{T}$ is from having Serre duality. Observe that in
the case of Example \ref{exm:3.6} the quotient category
$\mathcal{T}/{\mathcal{T}_r}$ coincides with the \emph{singularity
category} of the artin algebra $A$ (\cite{Or, CZ}); similarly, in
the (commutative) case of Theorem \ref{thm:3.10} the quotient
category $\mathcal{T}/{\mathcal{T}_r}$ coincides with the
singularity category of the corresponding (commutative) projective scheme.\\
(2) Note that there is an analogue of Reiten-Van den Bergh's theorem
for abelian categories in \cite{LZ}. Roughly speaking, a Hom-finite
abelian category has the so-called \emph{Auslander-Reiten duality}
if and only if it has enough Auslander-Reiten sequences. We expect
that there exists a reasonable notion of \emph{generalized
Auslander-Reiten duality}  for an arbitrary Hom-finite abelian
category. This generalized duality might be useful in the study of
an abelian category which does not have enough Auslander-Reiten
sequences. For example, the category of finite dimensional comodules
over a coalgebra is often such an abelian category.
\end{rem}

 \vskip 30pt

 \begin{appendix}

 \section{Generalized Serre duality on linear category}

Throughout this appendix, let $R$ be a commutative artinian ring.
Let $\mathcal{C}$  be an $R$-linear category which is Hom-finite. We
do not assume that $\mathcal{C}$ is an additive category. Use the
notation as in the introduction. For example, for two objects $X$
and $Y$, $(X, Y)$ stands for ${\rm Hom}_\mathcal{C}(X, Y)$; $D={\rm
Hom}_R(-, E)$ is the Matlis duality, where $E$ is the minimal
injective cogenerator of $R$.

\subsection{} In this subsection, we introduce the notion of
generalized Serre duality on a linear category. One might compare
with the treatment in \cite[Proposition 3.4]{BK}, \cite[p.301]{RV}
and \cite[section 3]{MS}.

Let $\mathcal{C}$ be a linear category as above. Recall the two full
subcategories $\mathcal{C}_r$ and $\mathcal{C}_l$  in the
introduction. We have the following basic result.

\begin{prop}\label{prop:A.1}
Use the notation as above. Then there is a unique functor (up to
isomorphism) $S \colon \mathcal{C}_r \rightarrow \mathcal{C}_l$ such
that for each  $X \in \mathcal{C}_r$ and $Y\in
 \mathcal{C}$, there is a natural $R$-linear isomorphism
 \begin{align*}
 \phi_{X, Y}\colon D(X, Y) \simeq (Y, S(X)).
\end{align*}
Moreover, the functor $S$ is an $R$-linear equivalence between
$\mathcal{C}_r$ and $\mathcal{C}_l$.
\end{prop}

\begin{proof} We will first construct the functor $S$. Take an object $X\in
\mathcal{C}_r$. Then there is an object $S(X)\in\mathcal{C}$ such
that there is an isomorphism of $R$-linear functors $\phi_X\colon
D(X, -) \simeq (-, S(X))$. We stress that this isomorphism is
required to preserve the $R$-linear structure.

Given any morphism $f\colon X \rightarrow Y$ in $\mathcal{C}_r$,
define the morphism $S(f)\colon S(X)\rightarrow S(Y)$ such that the
following diagram commutes
\[\xymatrix{
 D(X, -) \ar[d]_-{D(f, -)} \ar[r]^{\phi_X} & (-, S(X)) \ar[d]^-{(-, S(f))}\\
             D(Y, -) \ar[r]^{\phi_Y} & (-, S(Y))
}\] Here we apply Yoneda Lemma. Moreover, the morphism $S(f)$ is
uniquely determined by the above commutative diagram. It is direct
to check that $S\colon \mathcal{C}_r \rightarrow \mathcal{C}$ is a
$R$-linear functor. Note that for each object $X\in \mathcal{C}_r$,
we have an isomorphism of $R$-linear functors $D(-, S(X)) \simeq
(X,-)$, and then we have $S(X) \in \mathcal{C}_l$. In other words,
we have obtained an $R$-linear functor $S\colon \mathcal{C}_r
\rightarrow \mathcal{C}_l$. By the construction of $S$, we observe
that there is a natural $R$-linear isomorphism
$$ \phi_{X, Y}\colon D(X, Y) \stackrel{\simeq}{\longrightarrow} (Y, S(X))$$
for $X \in \mathcal{C}_r$ and $Y \in \mathcal{C}$. Here we put
$\phi_{X,Y}=\phi_X(Y)$. The naturalness of $\phi_{X, Y}$ in $Y$ is
obvious, while its naturalness in $X$ is a direct consequence of the
action of the functor $S$ on morphisms.

For the uniqueness of the functor $S$, assume that there is another
functor $S'\colon \mathcal{C}_r \rightarrow \mathcal{C}_l$ such that
there is a natural $R$-linear isomorphism $\phi'_{X, Y}\colon D(X,
Y) \simeq (Y, S'(X))$. In particular, there is an isomorphism of
$R$-linear functors $D(X, -)\simeq (-, S'(X))$. However by the above
proof, we have seen $\phi_X \colon D(X, -)\simeq (-, S(X))$.
Therefore $(-, S(X)) \simeq (-, S'(X))$, and hence Yoneda Lemma
implies that there is a unique isomorphism $S(X) \simeq S'(X)$. It
is routine to check that this gives rise to a  natural isomorphism
between the $R$-linear functors $S$ and $S'$.
\par \vskip 3pt

 We will show that the functor $S\colon \mathcal{C}_r\rightarrow
\mathcal{C}_l$ is fully faithful and dense,  and then  it is a
$R$-linear equivalence of categories. For this end, we need to
introduce some notation. For any $X\in \mathcal{C}_r$ and $Y \in
\mathcal{C}$, define a non-degenerated bilinear form
 $$(-,-)_{X, Y} \colon (X, Y) \times (Y, S(X)) \longrightarrow E$$
 such that $(f, \; g)_{X, Y}=\phi_{X, Y}^{-1}(g)(f)$. By the
 naturalness of $\phi_{X, Y}$, we have
 \begin{align}
 (f\circ \theta, \;g)_{X', Y}&=(f, \;S(\theta)\circ g)_{X, Y} \\
  (f,\; g\circ \gamma)_{X, Y'}&=(\gamma \circ f,
 \; g)_{X, Y},
\end{align}
where $\theta\colon X'\rightarrow X$ and $\gamma \colon
Y'\rightarrow Y$ are arbitrary morphisms. Note that in (A.1), $X,
X'\in \mathcal{C}_r$, $f\colon  X\rightarrow Y$ and $g\colon Y
\rightarrow S(X')$; in (A.2), $X\in \mathcal{C}_r$, $f\colon X
\rightarrow Y'$ and $g\colon Y \rightarrow S(X)$.

\vskip 3pt

We will call the bilinear form $(-, -)$ defined above  the
\emph{associated bilinear form} to  $S$. Note that the following
direct consequence of the above definition: $F(f)=(f, \phi_{X,
Y}(F))_{X, Y}$ for each $F \in D(X, Y)$ and $f\in (X, Y)$.\par

 Let $X, Y \in \mathcal{C}_r$. We
claim that  for each morphism $f\in (X, Y)$ and $g\in (Y, S(X))$ we
have
\begin{align}
(f, \; g)_{X, Y}=(g,\; S(f))_{Y,S(X)}.\end{align}

In fact, we have
\begin{align*}
(f, \; g)_{X,Y}&= (f \circ {\rm Id}_X,\; g)_{X, Y} \quad\quad  \mbox{use (A.1)}\\
&=({\rm Id}_X, \; S(f)\circ
g)_{X,X}  \quad \mbox{use (A.2)}\\
&= (g, \; S(f))_{Y,S(X)}.
\end{align*}
Consider the following composite of $R$-linear isomorphisms
\begin{align*}
\Phi \colon (X, Y) \stackrel{D(\phi_{X, Y})^{-1}}{\longrightarrow}
D(Y, S(X)) \stackrel{\phi_{Y, S(X)}}{\longrightarrow} (S(X), S(Y)).
\end{align*}
Note that  $D(\phi_{X, Y})^{-1} (f)(g)= (f, g)_{X, Y}$ and $F(g)=
(g, \phi_{Y, S(X)}(F))_{Y, S(X)}$ for all $F\in D(Y, S(X))$. Then it
is direct to see that $(f, g)_{X, Y}=(g, \Phi(f))_{Y, S(X)}$.
Compare this with (A.3), and note that the bilinear form is
non-degenerated. We infer that $S(f)=\Phi(f)$. Consequently the
functor $S$ is fully faithful.

To see that $S$ is dense, let $X' \in \mathcal{C}_l$. Note that the
functor $D(-, X')$ is representable. Assume that $D(-, X')\simeq (X,
-)$. Then $D(X,-) \simeq (-, X')$, and this implies that $X \in
\mathcal{C}_r$. Therefore we have
\begin{align*}
(-, S(X)) \stackrel{\phi_X^{-1}}\simeq D(X, -) \simeq (-, X').
\end{align*}
By Yoneda Lemma $S(X) \simeq X'$. Then the functor  $S$ is dense,
completing the proof. \end{proof}
\par \vskip 10pt

Following \cite{RV} we will introduce the notion of trace function.
For each object $X\in \mathcal{C}_r$, define its \emph{trace
function} to be an $R$-linear map
  $${\rm Tr}_X: (X, S(X))\longrightarrow E$$
such that ${\rm Tr}_X(f)=({\rm Id}_X, f)_{X, X}$. Then by (A.2), we
have $(f, g)_{X, Y}={\rm Tr}_X(g\circ f)$. By (A.3), we get
\begin{align}
{\rm Tr}_X(g\circ f)={\rm Tr}_Y(S(f)\circ g).
\end{align}

\vskip 10pt

 Let $\mathcal{C}$ be a Hom-finite $R$-linear category. We call the above obtained sextuple
$$\{S,\;  \mathcal{C}_r,\; \mathcal{C}_l, \; \phi,\; (-, -),\; {\rm
Tr}\}$$ the \emph{generalized Serre duality} on the category
$\mathcal{C}$. The functor $S$ is called the \emph{generalized Serre
functor} of $\mathcal{C}$, where $\mathcal{C}_r\subseteq
\mathcal{C}$ (\emph{resp}. $\mathcal{C}_l\subseteq \mathcal{C}$) is
referred as the \emph{domain}  (\emph{resp}. the \emph{range}) of
the generalized Serre functor.

We say that the category $\mathcal{C}$ \emph{has right Serre
duality} (\cite{RV}) provided that $\mathcal{C}_r=\mathcal{C}$. In
this case, we call $S$ the \emph{right Serre functor} of
$\mathcal{C}$. Similarly, we define the notion of \emph{having left
Serre duality}. The category $\mathcal{C}$ \emph{has Serre duality}
(\cite{BK}) provided that it has both right and left Serre duality.
This is equivalent to $\mathcal{C}_r=\mathcal{C}=\mathcal{C}_l$.

\subsection{}
 Let $\mathcal{C}$ be an $R$-linear category which is Hom-finite. Let $\mathcal{C}_r$
and $\mathcal{C}_l$ be the domain and the range of its generalized
Serre functor, respectively. In this subsection we will consider the
generalized Serre duality  on both $\mathcal{C}_r$ and
$\mathcal{C}_l$. \par \vskip 5pt

 Consider  the following two conditions on the category $\mathcal{C}$:
\begin{enumerate}
\item[(C)] for each $X, X' \in \mathcal{C}_r$ and $Z\in \mathcal{C}$,
there exists a morphism $s\colon Z' \rightarrow Z$ with $Z'\in
\mathcal{C}_r$ such that $(X, Z') \stackrel{(X, s)}{\simeq} (X, Z) $
and $(Z', X') \stackrel{(s, X')}{\simeq} (Z, X')$;
\item[(C')] for each $Y, Y'\in \mathcal{C}_l$ and $Z\in \mathcal{C}$,
there exists a morphism $s\colon Z \rightarrow Z'$ with $Z'\in
\mathcal{C}_l$ such that $(Y, Z) \stackrel{(Y, s)}{\simeq} (Y, Z')$
and $(Z', Y') \stackrel{(s, Y')}{\simeq} (Z, Y')$.
\end{enumerate}

We have the following result, which is inspired by \cite[Theorem
3.4]{H3}.

\begin{lem}\label{lem:A.2}
Let $\mathcal{C}$ be as above. Then we have \begin{enumerate}\item
if $\mathcal{C}_l\subseteq \mathcal{C}_r$, then
$(\mathcal{C}_r)_r=\mathcal{C}_r$;  \item conversely, if
$(\mathcal{C}_r)_r=\mathcal{C}_r$ and the category $\mathcal{C}$
satisfies the condition (C), then $\mathcal{C}_l \subseteq
\mathcal{C}_r$.\end{enumerate}
\end{lem}

\begin{proof} (1) is obvious from the definition.\par
 To
see (2), let $X \in \mathcal{C}_r$. Since
$(\mathcal{C}_r)_r=\mathcal{C}_r$, then we have an isomorphism of
$R$-linear functors on $\mathcal{C}_r$: $D(X, -) \simeq (-, X')$,
where $X'\in \mathcal{C}_r$. We claim that this isomorphism $D(X, -)
\simeq (-, X')$ can be extended on $\mathcal{C}$. Then we get
$S(X)=X'$, where $S$ is the generalized Serre functor on
$\mathcal{C}$. By Proposition \ref{prop:A.1} we have
$\mathcal{C}_l=S(\mathcal{C}_r)$. Then we conclude that
$\mathcal{C}_l \subseteq \mathcal{C}_r$.\par

In fact, for each $Z \in \mathcal{C}$, take $s\colon Z' \rightarrow
Z$ to be the morphism in the condition (C). Then we have a sequence
of isomorphisms
\begin{align*}
D(X, Z) \stackrel{D(X, s)}{\simeq} D(X, Z') \simeq (Z', X')
\stackrel{(s, X')}{\simeq} (Z, X').
\end{align*}
Note that this composite isomorphism is  given by the trace function
on $X$; more precisely, the inverse of this isomorphism is given by
$f \longmapsto (g \mapsto {\rm Tr}_X (f\circ g))$; here ${\rm Tr}$
is the trace function given by the generalized Serre duality on
$\mathcal{C}_r$. Then we get a natural isomorphism of $R$-linear
functors on $\mathcal{C}$: $D(X, -) \simeq (-, X')$. We are done.
\end{proof}
\vskip 5pt

 Dually, we have the following result.

\begin{lem}\label{lem:A.3}
Let $\mathcal{C}$ be as above. Then we have
\begin{enumerate}\item if $\mathcal{C}_r\subseteq \mathcal{C}_l$,
then $(\mathcal{C}_l)_l=\mathcal{C}_l$; \item conversely, if
$(\mathcal{C}_l)_l=\mathcal{C}_l$ and the category $\mathcal{C}$
satisfies the condition (C'), then $\mathcal{C}_r \subseteq
\mathcal{C}_l$.\hfill $\square$
\end{enumerate}
\end{lem}

\vskip 5pt

 The following result seems to be of interest.

\begin{prop}\label{prop:A.4}
Let $\mathcal{C}$ be a  Hom-finite $R$-linear category. Let
$\mathcal{C}_r$ and $\mathcal{C}_l$  the domain and the range of the
generalized Serre functor on $\mathcal{C}$, respectively. Assume
that the category $\mathcal{C}$ satisfies the conditions (C) and
(C'). Then the following statements are equivalent:
\begin{enumerate}
\item the category $\mathcal{C}_r$ has Serre duality;
\item  we have $\mathcal{C}_r=\mathcal{C}_l$;
\item the category $\mathcal{C}_l$ has Serre duality.
\end{enumerate}
\end{prop}

\begin{proof} The implications  (2) $\Rightarrow$ (1) and (1)
$\Rightarrow$ (3)  follow from Lemma \ref{lem:A.2}(1) and Lemma
\ref{lem:A.3}(1), respectively.

To see (1) $\Rightarrow$ (2), we apply Lemma \ref{lem:A.2}(2) and
then we get $\mathcal{C}_l\subseteq \mathcal{C}_r$. Therefore
$(\mathcal{C}_l)_l\supseteq \mathcal{C}_l\cap (\mathcal{C}_r)_l$.
Here we use an easy fact that for any full subcategory $\mathcal{D}$
of a category $\mathcal{C}$, one has $\mathcal{D}_l\supseteq
\mathcal{C}_l\cap \mathcal{D}$. However by (1),
$(\mathcal{C}_r)_l=\mathcal{C}_r$. Hence
$(\mathcal{C}_l)_l=\mathcal{C}_l$. Now applying Lemma
\ref{lem:A.3}(2), we get $\mathcal{C}_r\subseteq \mathcal{C}_l$.
Then we are done. Similarly we prove the implication (2)
$\Rightarrow$ (3).
\end{proof}
 \end{appendix}

\vskip 10pt

 \noindent {\bf Acknowledgements:} \quad  The author would like to
thank Prof. Michel Van den Bergh, Prof. Pu Zhang  and
Prof. Bin Zhu  and Dr. Adam-Christiaan van Roosmalen  
for their helpful comments.

\bibliography{}

\vskip 10pt

 {\footnotesize \noindent Xiao-Wu Chen, Department of
Mathematics, University of Science and Technology of
China, Hefei 230026, P. R. China \\
Homepage: http://mail.ustc.edu.cn/$^\sim$xwchen \\
\emph{Current
address}: Institut fuer Mathematik, Universitaet Paderborn, 33095,
Paderborn, Germany}

\end{document}